\documentstyle[secteqno]{article}      
\title{On the evaluation of the norm of an integral operator associated with the stability of 
 one-electron atoms}  
\author{V.I.Burenkov and W.D.Evans}      

\newtheorem{theorem}{Theorem}[section]
\newtheorem{lemma}[theorem]{Lemma}

\newtheorem{remark}[theorem]{Remark}
\newcommand{\beq}{\begin{equation}}
\newcommand{\bd}{\begin{definition}}
\newcommand{\bex}{\begin{example}}
\newcommand{\bc}{\begin{corollary}}
\newcommand{\bl}{\begin{lemma}}
\newcommand{\bt}{\begin{theorem}}
\newcommand{\br}{\begin{remark}}
\newcommand{\eeq}{\end{equation}}
\newcommand{\ed}{\end{definition}}
\newcommand{\eex}{\end{example}}
\newcommand{\ec}{\end{corollary}}
\newcommand{\el}{\end{lemma}}
\newcommand{\et}{\end{theorem}}
\newcommand{\er}{\end{remark}}

\begin{document}             

\maketitle                   

\begin{abstract}
The norm of an integral operator occurring in the partial wave decomposition of 
an operator $B$ introduced by Brown and Ravenhall in a model for relativistic one-electron
atoms is determined. The result implies that $B$ is non-negative and has no eigenvalue at 
$0$ when the nuclear charge does not exceed a specified critical value.
\end{abstract}

\section{Introduction}
\par
The operator referred to in the title is defined on $L^2(0,\infty)$ by
\beq
(T\phi)(x) := \int_0^{\infty} t(x,y)\phi(y)dy,\;\;\;0<x<\infty,
\eeq
where,
\begin{eqnarray}
t(x,y) =& \frac{1}{2}&\!\! \left\{\sqrt{\frac{\sqrt{x^2+1}+1}{x^2+1}}
g_0(x/y) \sqrt{\frac{\sqrt{y^2+1}+1}{y^2+1}} \right.\nonumber \\
& +  &\left. \sqrt{\frac{\sqrt{x^2+1}-1}{x^2+1}} g_1(x/y)
\sqrt{\frac{\sqrt{y^2+1}-1}{y^2+1}}\right\} 
\end{eqnarray}
with
\[
g_0(u)=\log \left|\frac{u+1}{u-1}\right|, g_1(u)=\frac{1}{2}
\left(u+\frac{1}{u}\right)\log\left|\frac{u+1}{u-1}\right| -1,\;\;\;u>0.
\]
To describe its role in relativistic stability, we require some background information.
It is well-known that the Dirac operator describing relativistic one-particle
systems is unbounded below, and that problems occur when it is extended as a
model for multi-particle systems. The root of the problem is that the Dirac
operator describes two different particles, namely electrons and positrons. In the
paper [2] Brown and Ravenhall overcame this difficulty by projecting onto the electron
subspaces only. Specifically, for a relativistic electron in the field of its
nucleus, their operator is
\beq
B:= \Lambda_+ (D_0 - \frac{e^2 Z}{|\cdot|}) \Lambda_+.
\eeq
The notation in (1.3) is as follows:
\begin{itemize}
\item
$D_0$ is the free Dirac operator
\[
D_0 = c \alpha \cdot \frac{\hbar}{i} \nabla + mc^2 \beta \equiv
\sum_{j=1}^3 c\frac{\hbar}{i} \alpha_j \frac{\partial}{\partial x_j} + mc^2 \beta
\]
where $ {\bf \alpha} := (\alpha_1,\alpha_2,\alpha_3) $ and $ \beta $ are the
Dirac matrices given by
\[
\alpha_j =\left( \begin{array}{cc}
0_2 & \sigma_{j}  \\
\sigma_j&0_2
\end{array}
\right) ,
\beta = \left( \begin{array}{cc}
1_2 & 0_2 \\
0_2 & -1_2
\end{array}
\right)
\]
with $ 0_2, 1_2 $ the zero and unit $ 2\times 2 $ matrices respectively, and
$\sigma_j$ the Pauli matrices
\[
\sigma_1=\left( \begin{array}{cc}
0&1\\
1&0
\end{array}
\right),
\sigma_2=\left( \begin{array}{cc}
0& -i\\
i&0
\end{array}
\right),
\sigma_3 = \left( \begin{array}{cc}
1&0\\
0& -1
\end{array}
\right);
\]
\item
$ \Lambda_+ $ denotes the projection of $ L^2(R^3) \otimes C^4 $ onto the
positive spectral subspace of $ D_0 $, that is $ \chi_{(0,\infty)} D_0 $.
If we set
\[
\hat{f}({\bf p}) \equiv [{\cal F}(f)]({\bf p}) :=
(\frac{1}{2\pi \hbar})^{3/2} \int_{R^3} e^{-i{\bf p}.{\bf x}/ \hbar}
f({ \bf x}) d{\bf x}
\]
for the Fourier transform of $ f $, then it follows that
\[
\hat{(\Lambda_+ f)}({\bf p}) = \Lambda_+({\bf p})\hat{f} ({\bf p})
\]
where
\beq
\Lambda_+({\bf p}) = \frac{1}{2} + \frac{c{\bf{ \alpha \cdot p}}+ mc^2 \beta}{2e(p)},
e(p) = \sqrt{c^2p^2+m^2c^4},
\eeq
with $ p = |{\bf p}| $;
\item
$ 2\pi \hbar $ is Planck's constant, $c$ the velocity of light, $m$ the
electron mass, $-e$ the electron charge and $Z$ the nuclear charge.
\end{itemize}
\par
The underlying Hilbert space in which $B$ acts is
\beq
{\cal H } = \Lambda_+ (L^2(R^3)\otimes C^4)
\eeq
and when it is bounded below,$ B $ generates a self-adjoint operator (also denoted by $B$)
which is the Friedrichs extension of the restriction
of $B$ to $ \Lambda_+(C_0^{\infty}(R^3)\otimes C^4) $.
\par
 The operator $B$ was later used by Bethe-Salpeter (see [1]) and is
referred to with their name in [3]. In [3] it is proved that $B$ is bounded
below if and only if the nuclear charge $Z$ does not exceed the critical value
\beq
Z_c= 2/[(\frac{\pi}{2} + \frac{2}{\pi})\alpha], \;\;\; \alpha=e^2/\hbar c,
\eeq
where $\alpha$ is Sommerfeld's fine structure constant; this range of $Z$
covers all natural elements. When $Z=Z_c$, it is proved in [3] that
\[
B \geq -(\frac{\pi^2 -4}{\pi^2 +4})mc^2.
\]
However, in [4] Hardekopf and Sucher had investigated $B$ numerically and
predicted that $B$ is in fact non-negative, and that, as for the Dirac operator, 
the ground state energy vanishes for $Z=Z_c$, i.e. $0$ is an eigenvalue of $B$. 
The first part of this prediction of Hardekopf and Sucher has recently been 
confirmed, but the second part contradicted, by Tix in [5]. 
Following the basic strategy in [3], but with a better choice of trial 
functions, Tix obtains a lower bound for $B$ which is shown to be 
positive for $Z \leq Z_c $, specifically
\beq
B \geq mc^2(1-\alpha Z - 0.002 \frac{Z}{Z_c}) > 0;
\eeq
the numerical factor is roughly 0.09 for $ Z=Z_c$.
\par
From the partial wave analysis of $B$, it is shown in [3] that for all
$\psi \in \Lambda_+(C_0^{\infty}(R^3) \otimes C^4) $ and with $(\cdot,\cdot)$ the standard inner product on
$ L^2(R^3) \otimes C^4$,
\begin{eqnarray}
(B\psi,\psi)& = &\sum_{(l,m,s)\in I}\left\{\int_0^{\infty}e(p)
|a_{l,m,s}(p)|^2dp \right.\nonumber \\
& - &\left. \frac{\alpha cZ}{\pi}\int_0^{\infty}\!\!\int_0^{\infty}
\overline{a}_{l,m,s}(p^{\prime})k_{l,s}(p^{\prime},p)
a_{l,m,s}(p)dpdp^{\prime}\right\}
\end{eqnarray}
where $I$ is the index set
\begin{eqnarray}
I=\{(l,m,s): l\in N_0,& m & = -l-1/2,\cdots, l+1/2, s=1/2,-1/2,\nonumber \\
& |m| & \neq l+1/2 \;\;\;\mbox{when}\,\, s=-1/2\},
\end{eqnarray}
the kernels $ k_{l,s}(p^{\prime},p)$ are given by
\beq
k_{l,s}(p^{\prime},p) = \frac{[e(p^{\prime})+e(0)]Q_l(\frac{1}{2}[\frac{p}{p^{\prime}}
+\frac{p^{\prime}}{p}])[e(p)+e(0)] + c^2p^{\prime}Q_{l+2s}(\frac{1}{2}[\frac{p}{p^{\prime}}
+\frac{p^{\prime}}{p}])p }{(2e(p)[e(p)+e(0)])^{1/2}(2e(p^{\prime})[e(p^{\prime})
+e(0)])^{1/2}}
\eeq
and
\beq
\sum_{(l,m,s)\in I} \int_0^{\infty}|a_{l,m,s}(p)|^2dp = \|\psi\|^2
:= \sum_{j=1}^4 \int_{R^3}|\psi_j|^2d{\bf x}.
\eeq
In (1.10) the $Q_l$ are the Legendre functions of the second kind. The strategy
in [3] was based on this decomposition of $B$ and the observation that
\beq
0 \leq k_{l,s}(p^{\prime},p) \leq k_{0,1/2}(p^{\prime},p), \;\;\;l\in N_0, s=1/2,-1/2.
\eeq
It would follow that $ B\geq 0 $ for $ Z \leq Z_c $ if and only if
\beq
\int_0^{\infty}\int_0^{\infty} a(p^{\prime})k_{0,1/2}(p^{\prime},p)
a(p) dp^{\prime}dp \leq \frac{\pi}{\alpha c Z_c} \int_0^{\infty} e(p)a(p)^2 dp
\eeq
for all non-negative measurable functions $a$. On setting
\beq
g_l(u) = Q_l(\frac{1}{2}[u+\frac{1}{u}]),\;\;\; l \in N_0
\eeq
\[
p=mcx, p^{\prime}=mcy, \phi(x) = \sqrt{e(mcx)} a(mcx),
\]
(1.13) becomes
\beq
\int_0^{\infty}\!\!\int_0^{\infty}t(x,y) \phi(x)\phi(y) dxdy \leq
(\frac{\pi^2}{4} +1) \int_0^{\infty} \phi(x)^2 dx,
\eeq
where $t(\cdot,\cdot)$ is defined in (1.2).
What we prove in this paper is that the constant $ \frac{\pi^2}{4}+1 $ in the
inequality (1.15) is sharp, and there are no extremal functions. Furthermore, we
show that these results imply that $ B\geq (1-\frac{Z}{Z_c})mc^2 $ for $ Z \leq Z_c $
and $ 0 $ is not an eigenvalue of $B$ when $ Z=Z_c $. Much of the analysis
continues to be valid for analogous inequalities defined by general kernel
functions $ t_{l,s} $ derived from the $k_{l,s} $.
\section{The main results}
\par
The operator $ T $ defined on $ L^2(0,\infty) $ by (1.1)is readily seen to be a 
bounded
symmetric operator and so
\beq
\sup \left\{ \frac{|(T\phi, \phi)|}{\|\phi\|^2}: \phi \in L^2(0,\infty), \phi \neq 0\right\}
= \|T| L^2(0,\infty) \rightarrow L^2(0,\infty) \|.
\eeq
Our main result is
\bt
Let $ T $ be defined by (1.1). Then
\begin{enumerate}
\item
$ \|T | L^2(0,\infty)\rightarrow L^2(0,\infty) \| = \frac{\pi^2}{4}+1$;
\item
the operator $T $ has no extremal functions.
\end{enumerate}
\et
\br
We recall that $ \phi $ is an extremal function of a bounded symmetric operator
$T$ if $ \phi \in L^2(0,\infty), \phi \neq 0\,\, a.e. $ and
$ \|T|L^2(0,\infty)\rightarrow L^2(0,\infty)\| = |(T\phi,\phi)|/\|\phi \|^2.$ Hence
Theorem 2.1 means the following: for all non-negative measurable functions $ \phi \in
L^2(0,\infty)$ with $ \phi \neq 0\,\, a.e. $
\beq
\int_0^{\infty}\!\!\int_0^{\infty} t(x,y) \phi(x) \phi(y)dxdy <
(\frac{\pi^2}{4}+1) \int_0^{\infty} \phi^2(x)dx
\eeq
and the constant $ \frac{\pi^2}{4}+1 $ is sharp. In turn, this implies
that the inequality (1.13) is valid, the constant $ \pi/\alpha c Z_c $ is sharp
and there is no function $ a\in L^2(0,\infty;e(p)dp) $ which is not null and for
which there is equality in (1.13). A consequence of Theorem 2.1 is
\er
\bt
Let $B$ be the self-adjoint operator generated in $\cal H $ by (1.3) and let
$Z_c $ be given by (1.6). Then
\begin{enumerate}
\item  if $ Z\leq Z_c, B \geq (1-\frac{Z}{Z_c})mc^2 $;
\item  if $Z=Z_c, 0 $ is not an eigenvalue of $B$;
\item  if $Z>Z_c, B $ is unbounded below.
\end{enumerate}
\et
{\bf Proof. }
Part 3 is proved in [3]. From (1.8),(1.10),(1.11) and (1.12) it follows that
\begin{eqnarray*}
(B\psi,\psi ) & \geq & \sum_{(l,m,s)\in I}\left\{
\int_0^{\infty}e(p)|a_{l,m,s}(p)|^2 dp \right.\\
& - & \left.\frac{\alpha cZ}{\pi} \int_0^{\infty}\int_0^{\infty}
|a_{l,m,s}(p^{\prime})| k_{0,1/2}(p^{\prime},p)
|a_{l,m,s}(p)| dpdp^{\prime}\right\} \\
& \geq & \sum_{(l,m,s) \in I} \int_0^{\infty}\left(1-\frac{Z}{Z_c}\right)e(p)
|a_{l,m,s}(p)|^2 dp \\
& \geq & \left(1-\frac{Z}{Z_c}\right)mc^2 \|\psi\|^2
\end{eqnarray*}
which establishes part 1. To prove part 2, suppose 0 is an eigenvalue of $B$ with
corresponding eigenfunction $ \psi $. By (1.12) and (1.13), all the summands on the
right-hand side of (1.8) are non-negative and consequently are zero as now $ B
\psi = 0 $. Also (1.11) implies that at least one of the functions $ a_{l,m,s},
a_{l_0,m_0,s_0} $ say, is not null. But this would imply that there is equality in
(1.13) with the function $ a=|a_{l_0,m_0,s_0}| $, contrary to Remark 2.2. Hence
the proof is complete.
\br
When the mass $m=0$, a proof of Theorem 2.3 is given in [3]. On setting $ p=x,
p^{\prime}=y, \phi(x)=\sqrt{cx} a(x) $ in (1.13) when $ m=0$ we obtain
\beq
\int_0^{\infty}\!\!\int_0^{\infty} t_0(x,y) \phi(x)\phi(y)dxdy \leq
\left(\frac{\pi^2}{4}+1\right) \int_0^{\infty} \phi(x)^2 dx
\eeq
where
\beq
t_0(x,y):= \frac{1}{2\sqrt{xy}}\left\{ g_0\left(\frac{x}{y}\right)
+ g_1\left(\frac{x}{y}\right) \right\}.
\eeq
We shall prove in Section 3 that the integral operator $ T_0 $ with kernel
$ t_0 $ satisfies Theorem 2.1, and thus yields the analogue of Theorem 2.3
in the case $ m=0 $.
\er
\br
Tix's lower bound (1.7) for $B$ is an improvement on that in Theorem 2.3(1).
If 0 is not in the essential spectrum $ \sigma_{ess}(B) $ of $B$ when 
$Z= Z_c $, as is the case when $Z<Z_c $ for $ \sigma_{ess}(B) = [mc^2,\infty) $
is established in [3, Theorem 2], then Parts 1 and 2 of Theorem 2.3 imply that
$B$ is strictly positive. However, no specific positive lower bound can be 
deduced from Theorem 2.1 alone.
\er
\section{Proof of Theorem 2.1}
\par
The starting point is the following simple result (cf[3, Section 2.3]). We shall denote by
$(\cdot,\cdot)$ and $\|\cdot\|$ the standard inner-product and
norm respectively in $L^2(0,\infty)$. It is sufficient to consider only
real-valued functions in $L^2(0,\infty)$ throughout this section.
\bl
Let $f,g,h $ be real-valued, measurable functions on $(0,\infty)$.
Moreover, let $g$ and $h$ be positive and
\beq
g(1/u) = g(u),\;\;\;0<u<\infty.
\eeq
Then
\beq
\int_0^{\infty}\!\!\int_0^{\infty}f(x)g\left(\frac{x}{y}\right)f(y)dxdy\leq
\int_0^{\infty}f(x)^2\left\{\int_0^{\infty}\frac{h(y)}{h(x)}g
\left(\frac{y}{x}\right)dy\right\}dx.
\eeq
Equality holds if and only if $ f(x) = A h(x) $ a.e. on $(0,\infty)$,
where $A$ is a constant.
\el
{\bf Proof.}
By the Cauchy-Schwarz inequality,
\begin{eqnarray*}
\lefteqn{\int_0^{\infty}\int_0^{\infty} f(x)g(x/y)f(y)dxdy} \\
& = & \int_0^{\infty}\!\!\int_0^{\infty}f(x)\sqrt{g(x/y)\frac{h(y)}{h(x)}}
f(y)\sqrt{g(y/x)\frac{h(x)}{h(y)}} dxdy \\
& \leq & \left(\int_0^{\infty}\!\!\int_0^{\infty} f(x)^2 g(x/y)
\frac{h(y)}{h(x)} dxdy\right)^{1/2}
\left(\int_0^{\infty}\!\!\int_0^{\infty} f(y)^2 g(y/x)
\frac{h(x)}{h(y)} dxdy\right)^{1/2}\\
& = & \int_0^{\infty}f(x)^2\left( \int_0^{\infty} g(y/x)
\frac{h(y)}{h(x)} dy\right)dx
\end{eqnarray*}
Equality holds if and only if, for some constants $\mu$ and $\lambda$
\[
\mu f(x) \sqrt{g(x/y)\frac{h(y)}{h(x)}} =
\lambda f(y) \sqrt{g(y/x)\frac{h(x)}{h(y)}}
\]
a.e. on $ (0,\infty) \times (0,\infty) $. This is equivalent to $ f(x) = Ah(x) $
a.e. on $(0,\infty)$, where $A$ is a constant.
\bl
Let $G$ be the symmetric operator defined on $L^2(0,\infty)$ by
\beq
Gf(x) := \int_0^{\infty} \frac{g(x/y)}{\sqrt{xy}} f(y)dy, \;\;\;0<x<\infty,
\eeq
where $g$ is a positive measurable function satisfying (3.1). Then
\beq
\|G|L^2(0,\infty) \rightarrow L^2(0,\infty)\| = \int_0^{\infty} g(u) \frac{du}{u}.
\eeq
Moreover, there are no extremal functions.
\el
{\bf Proof.}
By Lemma 3.1 with $h(u) = 1/u$ we get
\begin{eqnarray*}
|(Gf,f)| & = & |\int_0^{\infty} \int_0^{\infty} \frac{f(x)}{\sqrt{x}}g(x/y)
\frac{f(y)}{\sqrt{y}} dxdy| \nonumber \\
& \leq & \int_0^{\infty} \left(\frac{f(x)}{\sqrt{x}}\right)^2\left(\int_0^{\infty}
\frac{x}{y}g(y/x)dy\right)dx \nonumber \\
& = & \int_0^{\infty} g(u)\frac{du}{u} \int_0^{\infty}f(x)^2 dx.
\end{eqnarray*}
Hence,
\beq
\|G|L^2(0,\infty)\rightarrow L^2(0,\infty)\| \leq \int_0^{\infty}g(u)\frac{du}{u}.
\eeq
Furthermore, equality in (3.5) can hold if and only if $f(x) = A/\sqrt{x} $ a.e.
on $(0,\infty)$. Since $ A/\sqrt{x}\not\in L^2(0,\infty) $ unless $A=0$, it
follows that for all $f \in L^2(0,\infty), f \neq 0\,\, a.e.$,
\beq
|(Gf,f)| < \left(\int_0^{\infty} g(u)\frac{du}{u}\right) \|f\|^2.
\eeq
\par
In order to establish the inequality converse to (3.5) we take $ f_{\delta}(x)
= \frac{\chi_{(1,\delta)}(x)}{\sqrt{x}} $ as a test function
, where $ \chi_{(1,\delta)}$ denotes the characteristic function of $(1,\delta),
1<\delta < \infty $. By l'Hospital's Rule, we have as $ \delta \rightarrow \infty $,
\begin{eqnarray*}
\|G| L^2(0,\infty) \rightarrow L^2(0,\infty) \| & \geq &
\lim_{\delta \rightarrow \infty} \frac{(Gf_{\delta},f_{\delta})}{\|f_{\delta}\|^2} \\
& = & \lim_{\delta \rightarrow \infty}\left\{ (\ln \delta)^{-1}
\int_1^{\delta}\left( \int_1^{\delta} g(y/x)\frac{dy}{y} \right)\frac{dx}{x} \right\} \\
& = & \lim_{\delta \rightarrow \infty} \left\{ \int_1^{\delta} g(y/\delta )\frac{dy}{y} +
\int_1^{\delta} g(\delta /x) \frac{dx}{x}\right\} \\
& = & \lim_{\delta \rightarrow \infty} \int_{1/\delta}^{\delta} g(u)\frac{du}{u} \\
& = & \int _0^{\infty} g(u)\frac{du}{u}.
\end{eqnarray*}
The equality (3.4) follows from (3.5). From (3.6) it follows that there is no
extremal function.
\bl
Let $ T_0 $ be the symmetric operator in $ L^2(0,\infty) $ defined by
\[
T_0 f(x) := \int_0^{\infty} t_0 (x,y) f(y)dy,
\]
where $ t_0 $ is given by (2.4). Then
\begin{enumerate}
\item
$\|T_0 |L^2(0,\infty) \rightarrow L^2(0,\infty) \| = \frac{\pi^2}{4} + 1;$
\item
there are no extremal functions.
\end{enumerate}
\el
{\bf Proof.}
The results follow from Lemma 3.2 since
\begin{eqnarray}
\int_0^{\infty} g_0(u)\frac{du}{u} & = & 2\int_0^1 g_0(u) \frac{du}{u} \nonumber \\
& = & 2\int_0^1 \log \left|\frac{u+1}{u-1}\right| \frac{du}{u} \nonumber \\
& = & 4 \int_0^1 \left(\sum _{k=0}^{\infty} \frac{u^{2k}}{2k+1}\right) du\nonumber \\
& = & 4 \sum_{k=0}^{\infty} \frac{1}{(2k+1)^2} \nonumber \\
& = & \frac{\pi^2}{2}
\end{eqnarray}
and
\begin{eqnarray}
\int_0^{\infty} g_1(u) \frac{du}{u} & = & 2\int_0^1 g_1(u) \frac{du}{u} \nonumber \\
& = & 2\int_0^1\left(\frac{1}{2}\left[u+\frac{1}{u}\right] \log \left|\frac{u+1}{u-1}\right|-1\right) \frac{du}{u} \nonumber \\
& = & 2 \lim_{\epsilon \rightarrow 0+, \delta \rightarrow 1-}
\left[\frac{1}{2}\left(u-\frac{1}{u}\right)\log \left|\frac{u+1}{u-1}\right|\right]_{\epsilon}^{\delta}\nonumber \\
& = & 2.
\end{eqnarray}
\bl
The operator $ T $ defined in (1.1) satisfies
\[
\|T|L^2(0,\infty) \rightarrow L^2(0,\infty) \| \geq \frac{\pi^2}{4} +1.
\]
\el
{\bf Proof.}
As in Lemma 3.2, we take $ f_{\delta}(x) = \frac{\chi_{(1,\delta)}(x)}{\sqrt{x}},
1< \delta < \infty, $ as a test function. By l'Hospital's rule we obtain
\begin{eqnarray*}
\|T|L^2(0,\infty) \rightarrow L^2(0,\infty)\|
& \geq & \lim_{\delta \rightarrow \infty}\frac{ (Tf_{\delta}, f_{\delta})}
{\|f_{\delta}\|^2} \\
& = & \lim_{\delta \rightarrow \infty} \left\{ (\ln \delta )^{-1}
\int_1^{\delta}\left( \int_1^{\delta} t(x,y)\frac{dy}{\sqrt{y}}\right)
\frac{dx}{\sqrt{x}}\right\}\\
& = & \lim_{\delta \rightarrow \infty}\left\{ \sqrt{\delta}
\int_1^{\delta}t(\delta,y)\frac{dy}{\sqrt{y}}\right. \\
& + &\left. \sqrt{\delta} \int_1^{\delta} t(x,\delta )\frac{dx}{\sqrt{x}})\right\}\\
& = & 2 \lim_{\delta \rightarrow \infty}
\int_{1/\delta}^1\delta t(\delta ,\delta u)\frac{du}{\sqrt{u}}.
\end{eqnarray*}
It is readily seen from (1.2) that for $ 1<\delta <\infty $,
\[
\frac{\delta t(\delta, \delta u)}{\sqrt{u}} \leq \frac{g_0(u)+g_1(u)}{u}
\in L(0,1)
\]
and
\[
\lim_{\delta \rightarrow \infty }\frac{\delta t(\delta, \delta u)}{\sqrt{u}}
= \frac {g_0(u)+g_1(u)}{2u}.
\]
Hence, by the Dominated Convergence Theorem, (3.7) and (3.8), we have
\[
\|T|L^2(0,\infty) \rightarrow L^2(0,\infty) \| \geq
\int_0^1 \frac{g_0(u)+g_1(u)}{u}du = \frac{\pi^2}{4} +1.
\]
\bl
For all functions $h_0, h_1 $ which are positive and measurable on $(0,\infty) $
\beq
\|T|L^2(0,\infty) \rightarrow L^2(0,\infty) \|\leq A(h_0,h_1),
\eeq
where
\begin{eqnarray*}
A(h_0,h_1)& = &\frac{1}{2} \sup_{0<x<\infty}
\left(\frac{\sqrt{x^2+1}+1}{x^2+1}\int_0^{\infty}
\frac{h_0(y)}{h_0(x)}g_0(y/x)dy\right.\nonumber \\
& + &\left.\frac{\sqrt{x^2+1}-1}{x^2+1}\int_0^{\infty}
\frac{h_1(y)}{h_1(x)}g_1(y/x)dy\right).
\end{eqnarray*}
The operator $ T\;\; $has an extremal function $ \phi $ if and only if
\beq
\phi(x) = A_0 h_0(x) \sqrt{\frac{x^2+1}{\sqrt{x^2+1}+1}}
= A_1h_1(x)\sqrt{\frac{x^2+1}{\sqrt{x^2+1}-1}}
\eeq
and
\beq
\frac{\sqrt{x^2+1}+1}{x^2+1}\int_0^{\infty}\frac{h_0(y)}{h_0(x)}g_0(y/x)dy
+ \frac{\sqrt{x^2+1}-1}{x^2+1}\int_0^{\infty}\frac{h_1(y)}{h_1(x)}g_1(y/x)dy
= A_3,
\eeq
for some non-zero constants $ A_1,A_2, A_3.$
\el
{\bf Proof.}
By Lemma 3.1,
\begin{eqnarray*}
\lefteqn{\int_0^{\infty}\!\!\int_0^{\infty}t(x,y)\phi(x)\phi(y)dxdy}\\
& = &\frac{1}{2}\left\{\int_0^{\infty}\!\!\int_0^{\infty}
\sqrt{\frac{\sqrt{x^2+1}+1}{x^2+1}}\phi(x)g_0(x/y)
\sqrt{\frac{\sqrt{y^2+1}+1}{y^2+1}}\phi(y) dxdy \right.\\
& + &\left. \int_0^{\infty}\!\!\int_0^{\infty}
\sqrt{\frac{\sqrt{x^2+1}-1}{x^2+1}}\phi(x)g_1(x/y)
\sqrt{\frac{\sqrt{y^2+1}-1}{y^2+1}}\phi(y) dxdy \right\} \\
& \leq & \frac{1}{2} \int_0^{\infty} \left(\frac{\sqrt{x^2+1}+1}{x^2+1}
\int_0^{\infty}\frac{h_0(y)}{h_0(x)}g_0(y/x)dy\right. \\
& +&\left. \frac{\sqrt{x^2+1}-1}{x^2+1}\int_0^{\infty}
\frac{h_1(y)}{h_1(x)}g_1(y/x)dy\right) \phi(x)^2 dx\\
& \leq & A(h_0,h_1) \|\phi\|^2.
\end{eqnarray*}
Moreover, the first inequality becomes an equality if and only if
\[
\sqrt{\frac{\sqrt{x^2+1}+1}{x^2+1}} \phi(x) = A_0h_0(x),
\sqrt{\frac{\sqrt{x^2+1}-1}{x^2+1}} \phi(x) = A_1h_1(x)
\]
a.e. on $(0,\infty) $, for some constants $ A_0,A_1 $. The second inequality
becomes an equality if and only if $ \phi(x) = 0 $ a.e. on the set of all
$ x \in (0,\infty) $ for which
\[
\frac{\sqrt{x^2+1}+1}{x^2+1}\int_0^{\infty}\frac{h_0(y)}{h_0(x)}g_0(x/y)dy
+ \frac{\sqrt{x^2+1}-1}{x^2+1}\int_0^{\infty}\frac{h_1(y)}{h_1(x)}g_1(y/x)dy
< A(h_0,h_1).
\]
Since an extremal function $ \phi $ is not null, this inequality can only be
satisfied on a set of zero measure. Consequently, (3.11) holds a.e. for some
constant $A_3 $.
\br
We note that for all functions $ h_0,h_1 $ which are positive and
measurable on $(0,\infty) $
\begin{eqnarray*}
\lefteqn{\liminf_{x \rightarrow \infty }
\frac{1}{2}\left(\frac{\sqrt{x^2+1}+1}{x^2+1} \int_0^{\infty}
\frac{h_0(y)}{h_0(x)}g_0(y/x)dy\right.}\\
& &\left.+ \frac{\sqrt{x^2+1}-1}{x^2+1}\int_0^{\infty}
\frac{h_1(y)}{h_1(x)}g_1(y/x)dy\right) \\
& &\geq \frac{\pi^2}{4}+1.
\end{eqnarray*}
Indeed, let
\[
\hat{h}_j(\xi): = \liminf_{x \rightarrow \infty} \frac{h_j(\xi x)}{h_j(x)},
\;\;\;0< \xi <\infty, j=0,1,
\]
where the $ \liminf $ can be finite or infinite. Then
\begin{eqnarray*}
\hat{h}_j(1/\xi) & = & \liminf_{x \rightarrow \infty} \frac{h_j(x/\xi)}{h_j(x)}\\
& = & \liminf_{y \rightarrow \infty}\frac{h_j(y)}{h_j(\xi y)} \\
& = & \frac{1}{\hat{h}_j(\xi)}.
\end{eqnarray*}
By Fatou's Theorem
\begin{eqnarray*}
\lefteqn{
\liminf_{x \rightarrow \infty }
\frac{1}{2}\left(\frac{\sqrt{x^2+1}+1}{x^2+1}\int_0^{\infty}
\frac{h_0(y)}{h_0(x)}g_0(y/x)dy\right.}\\
& &\left.+ \frac{\sqrt{x^2+1}-1}{x^2+1}\int_0^{\infty}
\frac{h_1(y)}{h_1(x)}g_1(y/x)dy \right) \\
& & =
\liminf_{x \rightarrow \infty }
\frac{1}{2}\left(\frac{\sqrt{x^2+1}+1}{x^2+1}
\int_0^{\infty}x \frac{h_0(ux)}{h_0(x)}g_0(u)du \right.\\
&  &\left.+ \frac{\sqrt{x^2+1}-1}{x^2+1}
\int_0^{\infty}x\frac{h_1(ux)}{h_1(x)}g_1(u)du\right)\\
&  & \geq
\frac{1}{2}\left(\int_0^{\infty}\liminf_{x\rightarrow \infty}
\left[\frac{h_0(ux)}{h_0(x)}\right] g_0(u)du
+\int_0^{\infty}\liminf_{x\rightarrow \infty}
\left[\frac{h_1(ux)}{h_1(x)}\right] g_1(u)du\right)\\
&  & =
\frac{1}{2}\left(\int_0^{\infty} \hat{h}_0(u)g_0(u)du +
\int_0^{\infty} \hat{h}_1(u)g_1(u)du\right).
\end{eqnarray*}
Furthermore, on substituting $ u=v- \sqrt{v^2-1} $ when $0<u<1 $ and
$ u=v+ \sqrt{v^2-1} $ when $ u > 1 $ we have
\begin{eqnarray*}
\int_0^{\infty} \hat{h}_j(u)g_j(u)du &=& \int_0^{\infty} \hat{h}_j(u)
Q_j\left( \frac{1}{2} \left[u+\frac{1}{u}\right] \right)du \cr
&=& \int_1^{\infty}\{\hat{h}_j(v-\sqrt{v^2-1})(v-\sqrt{v^2-1})\\
&+& \hat{h}_j(v+\sqrt{v^2+1})(v+\sqrt{v^2+1})\}
Q_j(v) \frac{dv}{\sqrt{v^2-1}} \\
&=& \int_1^{\infty}\left\{ \frac{1}{\hat{h}_j(v+\sqrt{v^2-1})(v+\sqrt{v^2-1})}\right.\\
& + &\left. \hat{h}_j(v+\sqrt{v^2-1})(v+\sqrt{v^2-1})\right\}
Q_j(v)\frac{dv}{\sqrt{v^2-1}} \\
&\geq & 2\int_1^{\infty} Q_j(v)\frac{dv}{\sqrt{v^2-1}}\\
&=& \int_0^{\infty}Q_j(\frac{1}{2}[u+\frac{1}{u}])\frac{du}{u} \\
&=& \int_0^{\infty}g_j(u)\frac{du}{u}.
\end{eqnarray*}
This verifies the assertion. We also note that equality holds if and only
if $ \hat{h}_j(u) = 1/u $ a.e. on $(0,\infty) $. Thus to prove that
$ A(h_0,h_1) \leq \frac{\pi^2}{4}+1 $, and hence complete the proof of
Theorem 2.1, we must choose $h_0 $ and $h_1 $ in such a way that
$ \hat{h}_0(u) = \hat{h}_1(u) = 1/u $ a.e. on $(0,\infty) $.
\er
\bl
For all $ \phi \in L^2(0,\infty), \phi(x) \neq 0\,\, a.e., $ we have
\beq
\int_0^{\infty} \int_0^{\infty} t(x,y) \phi(x) \phi(y) dxdy
< C \int_0^{\infty} \phi(x)^2 dx,
\eeq
where
\beq
C = \sup_{0<x<\infty} F(x)
\eeq
and
\beq
F(x) = \frac{\pi}{2} (\sqrt{x^2+1}+1)\frac{\arctan x}{x} +
\frac{(\sqrt{x^2+1}-1)x}{x^2+1}.
\eeq
\el
{\bf Proof.}
We apply Lemma 3.5 with the choice (cf[5])
\beq
h_0(x)=\frac{x}{x^2+1},\;\;\; h_1(x)=\frac{1}{x}.
\eeq
>From (3.8),
\[
\int_0^{\infty}h_1(y)g_1(y/x)dy = 2.
\]
Also, on using Cauchy's Residue Theorem, we obtain
\begin{eqnarray*}
\int_0^{\infty}h_0(y)g_0(y/x)dy
&=& \int_0^{\infty}\frac{y}{y^2+1}\log \left|\frac{x+y}{x-y}\right| dy \nonumber \\
&=& \frac{1}{2}\int_{-\infty}^{\infty}\frac{y}{y^2+1}\log \left|\frac{x+y}{x-y}\right|dy \nonumber \\
&=& \frac{x^2}{2}\int_{-\infty}^{\infty}\frac{u}{(xu)^2+1}\log \left|\frac{1+u}{1-u} \right|du \nonumber \\
&=& \frac{x^2}{2}\Re \left[\int_{-\infty}^{\infty}\frac{u}{(xu)^2+1}\log \left(\frac{1+u}{1-u} \right) du\right]\nonumber \\
&=& \pi \Re \left[\frac{i}{2}\log \left(\frac{x+i}{x-i} \right)\right]\nonumber \\
&=& \pi \arctan x.
\end{eqnarray*}
Thus (3.13) is confirmed. Since the equality (3.10) is not satisfied by the
choice of $ h_0, h_1 $ in (3.15) for any constants $A_0,A_1$, it follows from Lemma 3.5 that there is
strict inequality in (3.12).
\par
The final link in the chain of arguments is
\bl
The constant $C$ in (3.13) is given by
\beq
C= \frac{\pi^2}{4}+1.
\eeq
\el
{\bf Proof.}
Since $ \lim_{x\rightarrow \infty} F(x) = \frac{\pi^2}{4}+1 $, we have that
$ C \geq \frac{\pi^2}{4}+1 $. To prove the reverse inequality we start by
substituting $ x = \tan 2v $ in $ F(x) $ to obtain
\[
F(\tan 2v) = \frac{\pi v + 4\sin^4\! v}{\tan v},\;\;\;0\leq v \leq \pi/4.
\]
We therefore need to prove that
\[
f(v) := \pi v + 4\sin^4\! v -(\frac{\pi^2}{4}+1)\tan v \leq 0, \;\;\;0\leq v\leq \pi/4.
\]
The following identities for the derivatives are easily verified :
\begin{eqnarray*}
f^{(1)}(v)& = &\pi + 16\sin^3\! v \cos v -(\frac{\pi^2}{4}+1)\sec^2\! v,\\
f^{(2)}(v)& = & 2\sin v \sec^3\! v g(v),
\end{eqnarray*}
where
\[
g(v) = 3\sin 2v + 3\sin 4v +\sin 6v -(\frac{\pi^2}{4}+1)
\]
and
\[
g^{(1)}(v) = 12 \cos 4v (1+\cos 2v).
\]
Since $ g(0)<0, g(\pi/8)>0, g(\pi/4)<0, g^{(1)}(v)>0 $ on $[0,\pi/8) $ and
$ g^{(1)}(v)< 0 $ on $(\pi/8,\pi/4] $ there exist $v_1,v_2 $ such that
$ 0<v_1 < v_2 < \pi/4, g(v_1) = g(v_2) = 0, g(v) < 0 $ on $ [0,v_1) $ and
$ (v_2, \pi/4] $, and $ g(v) > 0 $ on $ (v_1,v_2) $. Thus $ f(0) = 0, f^{(1)}(0)
< 0, f^{(2)}(0) = 0, f(\pi/4)=0, f^{(1)}(\pi/4) >0, $ and
$ f^{(2)}(\pi/4) < 0 $. Moreover, $ f^{(2)} $ vanishes at $v_1 $ and $ v_2 $,
is negative on $ (0,v_1) $ and $ (v_2,\pi) $, and positive on $ (v_1,v_2) $.
In particular, it follows that $ f^{(1)} $ is negative on $[0,v_1] $ and
positive on $ [v_2,\pi/4] $.
\par
Suppose that $ f(\xi) = 0 $ for some $ \xi \in (0,\pi/4) $. From $ f(0) =
f(\xi) = f(\pi/4) = 0 $, and the last sentence of the previous paragraph,
it follows that there exist $ \eta_1, \eta_2 $ such that $ f^{(1)}(\eta_1) =
f^{(1)}(\eta_2) = 0 $ and $ v_1< \eta_1 <\eta_2 <v_2 $. Consequently,
there exists $ v_3 \in (v_1,v_2) $ such that $ f^{(2)}(v_3) = 0 $ which
is contrary to what was established in the previous paragraph. Thus
$ f(v) \neq 0 $ on $ (0,\pi/4) $. Since $ f(0) = f(\pi/4) = 0 $ and
$ f^{(1)}(0)<0, f^{(1)}(\pi/4) > 0 $, it follows that $ f(v) \leq 0 $
on $ [0,\pi/4] $. The proof is therefore complete. \\
{\bf Proof of Theorem 2.1.}
Part 1 follows from Lemmas 3.4,3.7 and 3.8, and part 2 from (3.12).\\
{\sc Acknowledgements}
We are grateful to Christian Tix for sending us a preprint of [5]. We
also record our thanks to the European Union for support under
the TMR grant FMRX-CT 96-0001.\\
\begin{center}
{\bf REFERENCES}
\end{center}
\begin{enumerate}
\item
Hans A.Bethe and Edwin E.Salpeter. Quantum mechanics of one- and two-electron
atoms. In S.Flugge, editor, {\em Handbuch der Physik, XXXV}, pages 88-436.
Springer, Berlin, 1st edition 1957.
\item
G.E.Brown and D.G.Ravenhall. On the interaction of two electrons. {\em Proc.Roy.Soc.
London} A, 208 (A 1095): 552-559, September 1951
\item
William Desmond Evans, Peter Perry and Heinz Siedentop. The spectrum of
relativistic one-electron atoms according to Bethe and Salpeter.
{\em Commun. Math. Phys.}, 178(3), 733-746 (1996).
\item
G.Hardekopf and J.Sucher. Critical coupling constants for relativistic wave
equations and vacuum breakdown in quantum electrodynamics. {\em Phys. Rev. A}
31(4), 2020-2029 (1985).
\item
Christian Tix. Strict positivity of a relativistic Hamiltonian due to Brown
and Ravenhall. Preprint.
\end{enumerate}
\vspace{0.50in}
School of Mathematics,\\
University of Wales, Cardiff,\\
Senghennydd Road,\\
Cardiff CF2 4YH.

\end{document}